\documentclass[11pt]{article}
\usepackage{amsmath, epsfig, cite, amssymb, amsfonts,latexsym}
\newcommand{\qed}{{\hfill\rule{4pt}{7pt}}}

\newtheorem{thm}{Theorem}[section]

\newtheorem{defi}[thm]{Definition}
\newtheorem{lem}[thm]{Lemma}

\def\pf{\noindent {\it Proof.} }
\def\MX{\mathbb{X}}
\def\MY{\mathbb{Y}}
\def\Z{\mathbb{Z}}
\def\MB{\mathbb{B}}
\def\us{\underset{\geq}{\Omega}}
\def\pio{\pi_{\omega}}

\def\pf{\noindent {\it Proof.} }

\numberwithin{equation}{section}

\makeatletter \@addtoreset{equation}{section} \makeatother

\begin{document}

\begin{center}
{\bf \large Partition Analysis and Symmetrizing  Operators}
\end{center}

 \vspace{0.3cm}
 \begin{center}
   {\bf Amy M. Fu}\\
   Center for Combinatorics\\
   Nankai University, Tianjin 300071, P.R. China\\
   Email: fmu@eyou.com
  \end{center}
 \begin{center}
  {\bf Alain Lascoux}\\
  Nankai University, Tianjin 300071, P.R. China\\
  Email: Alain.Lascoux@univ-mlv.fr\\
  CNRS, IGM Universit\'e  de
  Marne-la-Vall\'ee\\
  77454 Marne-la-Vall\'ee Cedex, France\\
  \end{center}

  \vspace{1cm}

  \noindent {\bf Abstract.}

Using a symmetrizing operator,
we give a new expression for the Omega operator used by MacMahon
in Partition Analysis, and given a new life by Andrews and his coworkers.
Our result is stated in terms of Schur functions.

\hspace{1.5cm}

In his book "Combinatory Analysis", MacMahon introduced an Omega
operator.  Recently, Andrews et al \cite{Andrews1, Andrews2, Andrews3} further
developed the theory of Partition Analysis. We show in theorem
\ref{theo} that the Omega operator can be expressed by a
symmetrizing operator. As a consequence, we can formulate:
$$\us \, \lambda^k/\prod_{x \in \MX}(1-x\lambda)\prod_{y \in \MY}(1-\frac{y}{\lambda})$$
in terms of Schur functions of $\MX$ and $\MY$ (and therefore in terms of the elementary symmetric functions in
$\MX$ and $\MY$).

Recall the definitions of MacMahon's Omega operator
$\us$ and of the symmetrizing operator $\pio$.
\begin{defi}
$$\us
\sum_{s_1=-\infty}^{\infty}\cdots\sum_{s_r=-\infty}^{\infty}A_{s_1,\cdots,s_r}\lambda_1^{s_1}\cdots\lambda_r^{s_r}
:=\sum_{s_1=0}^{\infty}\cdots\sum_{s_r=0}^{\infty}A_{s_1,\cdots,s_r}.
$$
\end{defi}
By iteration, it is sufficient to treat the case of one variable $\lambda$ only .

\begin{defi} \cite{Lascoux}
Given $\MX=\{x_1,x_2,\cdots,x_n\}$ of cardinality $|\MX|=n$, the
symmetrizing operator $\pio$ is defined by:
$$
\forall f(x_1,\ldots,x_n),\ \pio f(x_1,\cdots,x_n)= \sum_{\sigma
\in \frak{S}(\MX)}\sigma
\left(\frac{f(x_1,\cdots,x_n)}{\triangle(\MX)}x_1^{n-1}\cdots
x_n^0\right),
$$
writing $\triangle(\MX)$ for the Vandermonde $\prod_{1\leq i <j \leq n}(x_i-x_j)$,
the sum being over all permutations $\sigma$ in the symmetric group $\frak{S}(\MX)$.
\end{defi}

Recall that complete symmetric functions $S^j(\MX)$ are defined by the generating function:
\begin{equation*}
\sum_{j=0}^{\infty}S^{j}(\MX)\lambda^{j}=\frac{1}{\prod_{i=1}^n(1-x_i\lambda)}.
\end{equation*}

 Complete symmetric functions are compatible with union of
 alphabets (denoted `+'). Given $\MY=\{y_1,y_2,\cdots,y_m\}$,
 we have:
  \begin{equation*}
  S^{n}(\MX+\MY)=\sum_{k=0}^nS^{k}(\MX)S^{n-k}(\MY).
  \end{equation*}

Schur functions have two classical expressions:
   \begin{eqnarray*}
   S_{\mu}(\MX)&=&
    \left|x_i^{\mu_j+j-1} \right|_{1 \leq i, j\leq n}/\triangle(\MX)
       =\left|S^{\mu_i-i+j}(\MX) \right|_{1 \leq i, j \leq n} ,
\end{eqnarray*}
where $\mu=[\mu_1,\cdots, \mu_n]$ with $\mu_1 \geq \mu_2\geq \cdots \geq \mu_n \geq 0$.

From the definition of $\pio$, we get \cite{Lascoux} :
 \begin{equation}\label{Schur}
 \pio\,  x_1^{\mu_1}\cdots x_n^{\mu_n}
    = \left|x_i^{\mu_j+j-1} \right|_{1 \leq i, j\leq n}/\triangle(\MX)
    =S_{\mu}(\MX).
 \end{equation}

This formula is still valid if $\mu\in \Z^n$,
$\mu_1>-n,\, \ldots, \mu_n>-1$~:
\begin{equation}\label{Schur2}
 \pio\,  x_1^{\mu_1}\cdots x_n^{\mu_n}   =S_{\mu}(\MX) ,
 \end{equation}
the Schur function $S_{\mu}$, still defined as the determinant
$|S^{\mu_i-i+j}|_{1 \leq i, j \leq n}$,
 being either null or equal to $\pm$ a Schur function indexed by a
 partition.

Symmetrizing first in $x_2,\ldots, x_n$, one also has, with the same hypotheses
on $\mu$~:
\begin{equation}\label{Schur2}
 \pio\,  x_1^{\mu_1}\, S_{\mu_2,\ldots,\mu_n}(x_2,\ldots, x_n)
           =S_{\mu}(\MX)\, .
 \end{equation}

\begin{lem}\label{Lemma}
Given $\MX$, $\MY$ and $k$ such that $0\leq k<|\MX|$, then one has:
\begin{equation}\label{le}
\pio \left(\sum_{j=0}^{\infty} x_1^{j-k}S^j(\MY) \right)
=\sum_{j=0}^{\infty} S^{j-k}(\MX)S^j(\MY)\, .
\end{equation}
\end{lem}

\pf  Since powers of $x_1$ range from $-k$ to $\infty$,
we can apply (\ref{Schur2}):
$$ \pio \left(\sum_{j=0}^{\infty} x_1^{j-k}S^j(\MY) \right) =
  \sum_{j=0}^{\infty}  S_{j-k,\, 0^{n-1}}(\MX)\, S^j(\MY) \, .$$
The terms such that $j<k$ are all null, being determinants with two identical rows,
and the sum reduces to the expression stated in the lemma.  \qed

Let us remark that the operator $\us$ relative to $x_1, \ldots, x_n$
can be obtained from the operator
$x_1, \ldots, x_{n+r}$, $r\geq 0$ by specializing $x_{n+1}, \ldots, x_{n+r}$ to $0$.
 Therefore we can suppose that $n$ be bigger than any given
 integer. This allows us in the following theorem to suppose that
 $n>k$.
\begin{thm} \label{theo}
 Given two alphabets
$\MX=\{x_1, x_2, \cdots, x_n\}$  and $\MY=\{y_1,y_2,\cdots, y_m\}$
of cardinality $n$ and $m$, let $\MB=1+\MY=\{1\}\cup \MY$. If
$k<n$, then we have:
\begin{eqnarray}  \label{and}
\pio\, \sum x_1^{j-k}S^j(\MB)
&=&\us\,  \frac{\lambda^k}{(1-x_1\lambda)\cdots(1-x_n\lambda)
    (1-{y_1 \over \lambda})\cdots (1-{y_m \over \lambda})}\nonumber\\[3pt]
 &=&\frac{\sum_{\mu}(-1)^{|\mu|}S_{\mu'}(\MB)S_{-k,\, \mu}(\MX)}
    { R(1,\MX\MB)},
\end{eqnarray}
where $R(1,\MX\MY)$ is equal to $\prod_{x \in \MX , y \in
\MY}(1-xy)$, and where the sum is over all partitions $\mu$ (the
sum is in fact finite). The vector $[-k,\,
\mu_1,\ldots,\mu_{n-1}]$ is denoted $-k,\, \mu$.
\end{thm}

\pf We first recall Cauchy's formula \cite[p. 65]{Macd}:
$$
R(1,\MX\MY)=\sum_{\mu}(-1)^{|\mu|}S_{\mu}(\MX)S_{\mu'}(\MY),
$$
where $\mu \rightarrow \mu'$ is the conjugation of partitions.

\begin{eqnarray*}
 \us\, \sum_{i,j=0}^{\infty}S^{i}(\MX)S^{j}(\MY)\lambda^{i-j+k}&=&
 \us\,  \frac{\lambda^k}{(1-x_1\lambda)\cdots(1-x_n\lambda)
(1-{y_1 \over \lambda})\cdots (1-{y_m \over \lambda})}\\
 &=&\sum_{i=0}^{\infty}S^{i}(\MX)\sum_{j=0}^{i+k}S^{j}(\MY)
 =\sum_{i=0}^{\infty}S^{i}(\MX)S^{i+k}(\MB)\\
 &=&\sum_{j=0}^{\infty}S^{j-k}(\MX)S^{j}(\MB)\, .
\end{eqnarray*}
On the other hand, lemma \ref{Lemma} allows us to write this last sum as
$\pio \left(\sum_{j=0}^{\infty}x_1^{j-k}S^j(\MB) \right)$. We shall now directly
compute the action of $\pio$ on $\sum_{j=0}^{\infty}x_1^{j-k}S^j(\MB)$,
denoting ${\MX} \setminus x_1=\{x_2,\ldots,x_n\}$.
\begin{eqnarray*}
   \pio\, \sum_{j=0}^{\infty}
   x_1^{j-k}S^j(\MB)&=&\pio\, x_1^{-k}\sum_{j=0}^{\infty}x_1^jS^j(\MB)\\
   &=&\pio\,\frac{x_1^{-k}}{R(1,x_1\MB)}
     =\pio\, \frac{x_1^{-k}R(1,(\MX \setminus x_1)\MB) }{R(1,\MX\MB)} \\[3pt]
   &=&\frac{\pio\, \left( x_1^{-k}
      \sum_{\mu}(-1)^{|\mu|}S_{\mu'}(\MB)\, S_{\mu}(\MX \setminus x_1) \right)}
       { R(1,\MX\MB)}\\[3pt]
   &=&  \frac{\sum_{\mu}(-1)^{|\mu|}S_{\mu'}(\MB)\,
                            S_{-k,\, \mu}(\MX)} { R(1,\MX\MB)}
\end{eqnarray*}
and the theorem is proved.  \qed

The result can be expressed in terms of elementary symmetric functions because
$
e_{i}(\MB)=e_{i}(\MY)+e_{i-1}(\MY)
$
and Schur functions are determinants in elementary symmetric functions.

In \cite[Theorem 2.1]{Andrews4}, the authors give a ``Fundamental
Recurrence'' for the numerator of (\ref{and}).

In \cite[Theorem 1.4]{Guo}, Guo-Niu Han expresses
the Omega operator in terms of Lagrange interpolation:

\begin{equation}\label{Han}
\us\, \frac{\lambda^k}{A(\lambda)B(\lambda)}
=\sum_{i=1}^n\frac{x_i^{n-1-k}}{(1-x_i)B(x_i)\prod_{j \neq i}(x_i-x_j)},
\end{equation}
where:
$$
A(\lambda)=\prod_{i=1}^n(1-x_i\lambda), B(\lambda)=\prod_{j=1}^m(1-y_j\lambda).
$$

Let us recall the definition  \cite{Lascoux} of
the Lagrange operator $L_{\MX}$:

\begin{defi}
\begin{equation*}
\forall f \in \frak{Sym}(1|n-1), \ \
L_{\MX}f(x_1,\ldots,x_n)=\sum_{x \in \MX}\frac{f(x,\MX \setminus
x)}{R(x,\MX \setminus x)},
\end{equation*}
where $\frak{Sym}(1|n-1)$ is the space of polynomials in
$x_1,x_2,\ldots,x_n$, symmetrical in $x_2,\ldots, x_n$,  and
$R(x,\MX \setminus x)=\prod_{x' \in \MX \setminus x}(x-x')$.
\end{defi}

We can express the Lagrange operator in terms of $\pio$.
\begin{lem}
$\forall f \in \frak{Sym}(1|n-1)$, we have:
\begin{equation}\label{oper}
\pio f(x_1,\ldots, x_n)=L_{\MX}\left(f(x_1,\ldots,x_n)\,
x_1^{n-1}\right).
\end{equation}
\end{lem}

\pf Elements of $f(x_1,x_2,\ldots,x_n)$ can be written as sums of
powers of $x_1$, with coefficients
 symmetrical in $x_1,\ldots,x_n$. Checking now that
 $$
 L_{\MX}(x_1^k \, x_1^{n-1})=\pio(x_1^k)=S^k(\MX),
 $$
 is immediate.  \qed

 Formula (\ref{oper}) shows that the Lagrange operator in formula
  (\ref{Han}) can be replaced by $\pio$, and therefore \cite[Theorem 1.4]{Guo}
 is a consequence of theorem \ref{theo}.

One does not need to suppose that all the $x_i$'s be distinct. Indeed,
in a Schur function, one may specialize $x_1,\ldots,x_k$ to the same value $a$.
This is more of a problem
in the Lagrange interpolation formula, where one has in that case to use
derivatives of different orders.

Let us finish with a small explicit example of the action of
$\pio$, for $n=2, \, m=1, \, k=1$.
 \begin{eqnarray*}
    \pio \left(\sum_{j=0}^{\infty} x_1^{j-1}S^j(\MB)\right)
    &=&\frac{\sum_{\mu}(-1)^{|\mu|}S_{\mu'}(\MB)S_{-1,\, \mu}(\MX)}
            { R(1,\MX\MB)}\\
            &=&\frac{-S_{1}(\MB)S_{-1,\, 1}(\MX)+S_{1,1}(\MB)S_{-1,\, 2}(\MX)}
                                                                { R(1,\MX\MB)}\\
            &=&\frac{(1+y)-y(x_1+x_2)}{
            (1-x_1)(1-x_2)(1-x_1y)(1-x_2y)}\\
            &=& \us\, \frac{\lambda}{(1-\lambda x_1)(1-\lambda x_2)(1-y/\lambda)}.
           \end{eqnarray*}

\vskip 2mm \noindent{\bf Acknowledgments.}

This work was done under the auspices of the National Science
Foundation of China. The second author
is partially supported by the EC's IHRP
Program, within the Research Training Network ``Algebraic Combinatorics
in Europe'', grant HPRN-CT-2001-00272.

 \end{document}